\documentclass[12pt]{article}

\usepackage[margin=1.20in]{geometry}
\usepackage{sectsty}
\usepackage{xcolor}
\usepackage{amssymb,amsmath}
\usepackage{amsthm}
\usepackage{hyperref}
\usepackage{mathrsfs}
\usepackage{textcomp}
\usepackage{enumitem}
\usepackage[abbrev,nobysame,non-compressed-cites]{amsrefs}
\usepackage{todonotes}

\hypersetup{
  colorlinks,
  citecolor=black,
  filecolor=black,
  linkcolor=black,
  urlcolor=black
}

\newtheorem{Theorem}{Theorem}[section]
\newtheorem{Lemma}[Theorem]{Lemma}
\newtheorem{Proposition}[Theorem]{Proposition}

\newtheorem{Remark}[Theorem]{Remark}

\newcommand{\R}{\mathbb{R}}

\begin{document}

\title{Quantitative homogenization of $k$-Hessian equations in periodic media}

\author{
Tianling Jin\footnote{T. Jin was partially supported by NSFC grant
12122120, and Hong Kong RGC grants GRF 16304125, GRF 16306320 and
GRF 16303822.}, \quad
YanYan Li\footnote{Y.Y. Li is partially supported by NSF grants
DMS-2247410 and DMS-2554267.}, \quad
Hung V. Tran\footnote{H.V. Tran is partially supported by NSF grant
DMS-2348305.}, \quad
Xushan Tu
}

\date{\today}

\maketitle

\begin{center}
\vspace{-1cm}
\emph{\small{Dedicated to Professor Yaotian Shen on the occasion of his 90th birthday.\medskip}}
\end{center}

\begin{abstract}
Let $1\leq k\leq n$, and let $u^\varepsilon$ and $u$ be the $k$-admissible solutions of
\[
\sigma_k(D^2u^\varepsilon)=f(x,x/\varepsilon),\qquad\sigma_k(D^2u)=\overline f(x)
\]
on a bounded convex domain $\Omega$, with the same Dirichlet data. Here, $f$ is uniformly positive and periodic in its second variable, and $\overline f(x)=\int_{\mathbb T^n}f(x,y)\,dy$. For $m=0,1$ and $0<\alpha\leq1$, we prove
\[
\|u^\varepsilon-u\|_{L^\infty(\Omega)}\leq C\varepsilon^{m+\alpha},\qquad \forall\,0<\varepsilon\leq1,
\]
provided that $u\in C^{m+2,\alpha}(\overline\Omega)$ and $f\in C_x^{m,\alpha}(\overline\Omega;C_y^{2,\gamma}(\mathbb T^n))$ for some $0<\gamma<1$. Periodic cell estimates give smooth dependence of the correctors on the background Hessian and the normalized density. The divergence-free property of the Newton tensor cancels the first macroscopic correction and yields the higher rate. Both exponents are optimal in general, already for $n=k=1$.


\end{abstract}

\section{Introduction}

In this paper, we study quantitative homogenization for the $k$-Hessian equation with positive periodic densities on a bounded convex domain $\Omega\subset\R^n$. Let $1\leq k\leq n$, $\mathbb T^n:=\mathbb R^n/\mathbb Z^n$ be the flat $n$-dimensional torus, and $\varphi\in C(\partial\Omega)$. Let $\mathcal S^{n\times n}$ be the set of symmetric matrices of size $n$. For $M \in \mathcal S^{n\times n}$ with eigenvalues $\lambda_1,\ldots,\lambda_n$, set
\[
\sigma_k(M):=\sum_{1\leq i_1<\cdots<i_k\leq n}\lambda_{i_1}\cdots\lambda_{i_k},\qquad\Gamma_k:=\{M\in\mathcal S^{n\times n}:\sigma_j(M)>0,\ 1\leq j\leq k\}.
\]
We work with $k$-admissible solutions. For a $C^2$ function $u$, this means $D^2u\in\overline\Gamma_k$. A classical admissible solution with positive right-hand side satisfies $D^2u\in\Gamma_k$.

Let $f:\overline\Omega\times\mathbb T^n\to(0,\infty)$ be continuous and periodic in its second variable, and set
\[
\overline f(x):=\int_{\mathbb T^n}f(x,y)\,dy.
\]
Let $u^\varepsilon$ be a continuous $k$-admissible viscosity solution of
\begin{equation}
\label{eq:oscillatory-problem}
\sigma_k(D^2u^\varepsilon)=f(x,x/\varepsilon)\quad\text{in }\Omega,\qquad u^\varepsilon=\varphi\quad\text{on }\partial\Omega.
\end{equation}
The effective problem is
\begin{equation}
\label{eq:effective-problem}
\sigma_k(D^2u)=\overline f(x)\quad\text{in }\Omega,\qquad u=\varphi\quad\text{on }\partial\Omega.
\end{equation}
The periodic cell-average identity in Section~\ref{subsec:periodic-solutions} identifies this effective equation. We quantify the convergence of $u^\varepsilon$ to $u$ in terms of the regularity of $u$ and the density in $x$. We obtain the rates $O(\varepsilon^\alpha)$ and $O(\varepsilon^{1+\alpha})$ under the respective regularity assumptions stated below.

The literature on homogenization is vast. Here we briefly recall some results on periodic homogenization of linear equations in nondivergence form and fully nonlinear uniformly elliptic equations. For periodic linear elliptic equations in nondivergence form, Avellaneda and Lin \cite{avellaneda1989compactness} established regularity estimates uniform in the microscopic period. Evans \cite{evans1992periodic} proved periodic homogenization for fully nonlinear elliptic equations using the perturbed test function method. Caffarelli and Souganidis \cite{caffarelli2010souganidis} obtained algebraic convergence rates for periodic homogenization of fully nonlinear uniformly elliptic equations. Under higher regularity assumptions, Kim and Lee \cite{kim2016higher} constructed higher-order expansions for periodic linear and fully nonlinear uniformly elliptic equations using interior and boundary layer correctors. For general linear uniformly elliptic equations in nondivergence form, Guo, Tran, and Yu \cite{GTY2020} proved that $O(\varepsilon)$ is the optimal rate. When the diffusion matrix $A$ has only one degree of freedom, that is $A(y)=C+a(y)M$ where $C, M$ are constant symmetric matrices and $a$ is a scalar periodic function, the rate is improved to $O(\varepsilon^2)$ by Guo, Sprekeler, and Tran \cite{GST2025}.

This uniformly elliptic framework does not apply directly to the $k$-Hessian equation for $k\geq2$, which is not uniformly elliptic on the whole admissible cone. For the Monge--Amp\`ere equation, corresponding to $k=n$, periodic solutions and their regularity were studied by Li \cite{li1990existence} and Caffarelli and Li \cite{caffarelli2004liouville}.  Moreover, Caffarelli and Li \cite{caffarelli2004liouville} proved an algebraic homogenization estimate. Their method can also be used to extend this estimate to normalized nonnegative locally finite periodic Borel measures; see Theorem~2.3 in \cite{jin2025liouville}. In our previous paper \cite{jin2026quantitative}, we established the convergence rates $O(\varepsilon)$ and $O(\varepsilon^{2})$ for the Monge--Amp\`ere equation under the corresponding regularity assumptions in $x$, and also the $O(\varepsilon)$ rate for locally weighted periodic Borel measures, with constants independent of the distribution of the microscopic measure. Independently, Cui, Gu, Le, Wang, and Zhou \cite{cui2026global} obtained an $O(\varepsilon)$ convergence rate for the Monge--Amp\`ere equation with bounded, uniformly positive periodic densities $f(x/\varepsilon)$ on smooth uniformly convex domains with $C^3$ boundary data.

\medskip 

Our main result is the following.

\begin{Theorem}
\label{thm:main-result}
Let $1\leq k\leq n$, $m\in\{0,1\}$, $\alpha\in(0,1]$, $\gamma\in(0,1)$, and $\lambda, \Lambda>0$. Assume that $f:\overline\Omega\times\mathbb T^n\to(0,\infty)$ satisfies $0<\lambda\leq f\leq\Lambda$. Let $u^\varepsilon$ and $u$ be the $k$-admissible solutions of \eqref{eq:oscillatory-problem} and \eqref{eq:effective-problem}, respectively. If $f\in C_x^{m,\alpha}(\overline\Omega;C_y^{2,\gamma}(\mathbb T^n))$ and $u\in C^{m+2,\alpha}(\overline\Omega)$, then
\begin{equation}
\label{eq:main-rate}
\|u^\varepsilon-u\|_{L^\infty(\Omega)}\leq C\varepsilon^{m+\alpha},\qquad \forall\,0<\varepsilon\leq1.
\end{equation}
The constant $C$ depends only on $n$, $k$, $\alpha$, $\gamma$, $\lambda$, $\Lambda$, the inradius of $\Omega$, $\operatorname{diam}\Omega$, and the indicated norms of $f$ and $u$.
\end{Theorem}

\begin{Remark}
The $k$-admissibility of $u$, the lower bound $\sigma_k(D^2u)=\overline f\geq\lambda$, and the upper bound for $|D^2u|$ place $D^2u(\overline\Omega)$ in a fixed compact subset of $\Gamma_k$. Indeed, the Newton--Maclaurin inequalities give a positive lower bound for every $\sigma_j(D^2u)$, $1\leq j\leq k$. For $k=n$, this is the uniform convexity used in the Monge--Amp\`ere case.
\end{Remark}

\begin{Remark}
The quadratic endpoint is sharp in the one-dimensional case $n=k=1$: for $\Omega=(0,1)$, $\varphi=0$, $f(x,y)=1+\frac12\cos(2\pi y)$, and $\varepsilon=N^{-1}$ with $N\in\mathbb N$, direct integration yields
\[
\|u^\varepsilon-u\|_{L^\infty((0,1))}=\frac{1}{4\pi^2}\varepsilon^2.
\]
The two exponents cannot be increased for the full class in Theorem~\ref{thm:main-result}, already for $n=k=1$; see \cite{jin2026quantitative}*{Example~3.3}.
\end{Remark}

\begin{Remark}
The assumption $u\in C^{3,1}(\overline\Omega)$ in Theorem \ref{thm:main-result} with $m=\alpha=1$ holds, for example, if $\Omega$ is uniformly convex, $\partial\Omega\in C^{4,\gamma}$, $\varphi\in C^{4,\gamma}(\partial\Omega)$, and $\overline f\in C^{2,\gamma}(\overline\Omega)$. Indeed, since $\overline f\geq\lambda>0$, the global estimates of Caffarelli, Nirenberg, and Spruck \cite{caffarelli1985dirichlet} for the $k$-Hessian Dirichlet problem and Schauder estimates give
\[
u\in C^{4,\gamma}(\overline\Omega)\subset C^{3,1}(\overline\Omega).
\]
\end{Remark}

The proof of Theorem \ref{thm:main-result} is based on the method developed in our previous paper \cite{jin2026quantitative}. We use the periodic cell problems associated with the background Hessian $A=D^2u$ and the normalized density $h=f/\overline f$. We mollify these parameters at scale $\varepsilon$ and use the corresponding correctors to construct global approximate solutions. At the next order, the divergence-free property of the Newton tensor ensures solvability of the cell problem for the third-order corrector without a macroscopic correction. The comparison principle then gives the two uniform convergence rates from the residual estimates.

The rest of the paper is organized as follows. Section~\ref{sec:periodic-cell-problems} introduces the periodic cell problems and records a comparison principle for approximate solutions. Section~\ref{sec:corrector-estimates} develops the corrector expansions and proves Theorem~\ref{thm:main-result}.

\section{Preliminaries}
\label{sec:periodic-cell-problems}

Throughout the paper, $c$ and $C$ denote positive constants whose values may change from line to line. We use the Einstein summation convention and write $M:N=M_{ij}N_{ij}$ for symmetric matrices. The oscillatory equation is understood in the $k$-admissible viscosity sense. The effective solution has the regularity stated in Theorem~\ref{thm:main-result}. We use the comparison principle in this admissible class. In particular, a smooth strictly admissible function with larger $k$-Hessian and smaller boundary values lies below a solution, and the reversed inequalities give an upper comparison function. To see this, first make the Hessian inequality strict by a small quadratic perturbation and then apply the viscosity touching inequalities.

\subsection{Periodic solutions}
\label{subsec:periodic-solutions}

For $A\in\Gamma_k$, write $q_A(y)=\frac12y\cdot Ay$. A periodic cell potential has the form $w=q_A+v$, where $v$ is $\mathbb Z^n$-periodic and $D^2w\in\Gamma_k$. The linearization of $\sigma_k$ is the Newton tensor
\[
T_{k-1}(M):=D\sigma_k(M).
\]
It is positive definite for $M\in\Gamma_k$ and satisfies
\begin{equation}
\label{eq:newton-divergence}
\partial_{y_j}T_{k-1}^{ij}(D_y^2w)=0\qquad\forall\,i=1,\ldots,n.
\end{equation}
This identity holds for smooth Hessian fields and, by approximation, in distributions for $C^2$ potentials; see the cell-average argument in \cite{li2026quadratic}.

For every periodic $\zeta\in C^2(\mathbb T^n)$, the same identity gives
\begin{equation}
\label{eq:periodic-cell-compatibility}
\int_{\mathbb T^n}\sigma_k(A+D_y^2\zeta)\,dy=\sigma_k(A).
\end{equation}
Indeed, by replacing $\zeta$ by $t\zeta$ in the left-hand side, differentiating in $t$, and integrating by parts, one obtains that the derivative vanishes, which proves the identity. Since $\sigma_k$ is a polynomial, this identity does not require admissibility. Differentiating \eqref{eq:periodic-cell-compatibility} in $A$, with $\zeta$ fixed, also yields
\begin{equation}
\label{eq:periodic-newton-average}
\int_{\mathbb T^n}T_{k-1}(A+D_y^2\zeta)\,dy=T_{k-1}(A).
\end{equation}

For a positive periodic density $h$ with $\int_{\mathbb T^n}h=1$, let $v_{A,h}$ denote the normalized solution of
\begin{equation}
\label{eq:parametric-corrector}
\sigma_k(A+D_y^2v_{A,h})=\sigma_k(A)h,\qquad\int_{\mathbb T^n}v_{A,h}\,dy=0,
\end{equation}
and write $w_{A,h}=q_A+v_{A,h}$. Theorem~\ref{thm:smooth-periodic-cell-problem} gives the existence and regularity needed below. In particular, \eqref{eq:periodic-cell-compatibility} identifies $\overline f$ as the effective density when $A=D^2u$ and $h=f/\overline f$.

\subsection{A comparison principle}

We shall use the following consequence of the comparison principle.

\begin{Proposition}
\label{prop:approximate-solution-principle}
Fix $\beta\in(0,2]$, $C_0>0$, and a compact set $\mathcal K\Subset\Gamma_k$. There exist $\varepsilon_0>0$ and $C>0$, depending only on $n$, $k$, $\mathcal K$, $C_0$, $\beta$, and $\operatorname{diam}\Omega$, such that the following holds for every $0<\varepsilon<\varepsilon_0$. Let $u^\varepsilon$ and $u$ be the solutions of \eqref{eq:oscillatory-problem} and \eqref{eq:effective-problem}, respectively. Suppose there exists $V^\varepsilon\in C^2(\Omega)\cap C(\overline\Omega)$ such that
\[
\|V^\varepsilon-u\|_{L^\infty(\overline\Omega)}\leq C_0\varepsilon^\beta,\qquad D^2V^\varepsilon(x)\in\mathcal K\qquad\forall\,x\in\Omega,
\]
and
\[
\left|\sigma_k(D^2V^\varepsilon)-f(x,x/\varepsilon)\right|\leq C_0\varepsilon^\beta\quad\text{in }\Omega.
\]
Then
\[
\|u^\varepsilon-u\|_{L^\infty(\Omega)}\leq C\varepsilon^\beta.
\]
\end{Proposition}

\begin{proof}
After a translation, we may assume that $\Omega$ lies in a ball centered at the origin with radius at most $2\operatorname{diam}\Omega$. Set
\[
W^\pm=V^\varepsilon\pm c\varepsilon^\beta|x|^2.
\]
On a fixed neighbourhood of $\mathcal K$ in $\Gamma_k$, the Newton tensor is uniformly positive definite. Since $D^2W^\pm=D^2V^\varepsilon\pm2c\varepsilon^\beta I$, we may choose $c$ sufficiently large and then $\varepsilon_0>0$ sufficiently small so that, for every $0<\varepsilon<\varepsilon_0$, $W^\pm$ are $k$-admissible and
\[
\sigma_k(D^2W^+)\geq f(x,x/\varepsilon)\geq\sigma_k(D^2W^-)\quad\text{in }\Omega.
\]
On $\partial\Omega$, where $u^\varepsilon=u$, the bound $\|W^\pm-u\|_{L^\infty(\overline\Omega)}\leq C\varepsilon^\beta$ gives
\[
W^+-C\varepsilon^\beta\leq u^\varepsilon\leq W^-+C\varepsilon^\beta.
\]
The comparison principle extends these inequalities to $\Omega$. Combining them with $\|W^\pm-u\|_{L^\infty(\overline\Omega)}\leq C\varepsilon^\beta$ proves the proposition.
\end{proof}

\section{Corrector expansions for positive densities}
\label{sec:corrector-estimates}

In this section, we adapt the corrector construction of \cite{jin2026quantitative} to prove Theorem~\ref{thm:main-result}. To motivate the construction, assume that the data are smooth and consider the formal expansion (ansatz)
\[
u^\varepsilon(x)=u(x)+\varepsilon a(x)+\varepsilon^2v(x,y)+\varepsilon^3\psi(x,y)+\cdots,
\]
where $y=x/\varepsilon$ and $v,\psi$ are periodic in $y$. Substituting this expansion into \eqref{eq:oscillatory-problem} and matching the zeroth-order terms, we obtain
\[
\sigma_k\bigl(D^2u(x)+D_y^2v(x,y)\bigr)=f(x,y).
\]
Since $\sigma_k(D^2u)=\overline f$, we take the normalized solution $v(x,y)=v_{A(x),h(x,\cdot)}(y)$, where $A(x)=D^2u(x)$ and $h(x,y)=f(x,y)/\overline f(x)$. Collecting the order-$\varepsilon$ terms gives
\[
W^{ij}\bigl(a_{ij}+v_{x_i y_j}+v_{x_j y_i}+\psi_{y_i y_j}\bigr)=0,
\]
where $W=T_{k-1}(A+D_y^2v)$. The terms containing $y$-derivatives have zero integral over $\mathbb T^n$ by integration by parts, periodicity, and \eqref{eq:newton-divergence}. It follows from \eqref{eq:periodic-newton-average} that
\[
T_{k-1}^{ij}(D^2u)a_{ij}=0\quad\text{in }\Omega.
\]
Matching the boundary values at order $\varepsilon$ gives $a=0$ on $\partial\Omega$, so the maximum principle yields $a\equiv0$. We therefore determine $\psi(x,\cdot)$ from
\[
W^{ij}\psi_{y_i y_j}=-2W^{ij}v_{x_i y_j}\quad\text{in }\mathbb T^n,\qquad\int_{\mathbb T^n}\psi(x,y)\,dy=0,
\]
where the normalization fixes the additive constant. The derivative $v_{x_i}$ contains contributions from both $A$ and $h$, and the same integration by parts cancels both.

For general linear elliptic equations in nondivergence form, the first macroscopic correction need not vanish, and the $O(\varepsilon)$ rate is generally optimal; see \cites{GTY2020,GST2025}.
Here the divergence-free property of $T_{k-1}(A+D_y^2v)$ forces $a\equiv0$.
In our setting, the Newton tensor replaces the cofactor matrix in the Monge--Amp\`ere argument in \cite{jin2026quantitative}.
 This cancellation yields the rate $O(\varepsilon^{1+\alpha})$ for $0<\alpha\leq1$ under the regularity assumptions in Theorem~\ref{thm:main-result}.

Li and Liang \cite{li2026quadratic} established existence and uniform estimates for the periodic $k$-Hessian cell problem with positive definite background matrices, building on the a priori estimates of Li \cite{li1990existence} for Hessian equations. We extend these estimates to background matrices in compact subsets of $\Gamma_k$. The estimates place the cell Hessians in a fixed compact subset of $\Gamma_k$ and allow us to differentiate the correctors with respect to their parameters without assuming convexity of $u$. 

\begin{Theorem}
\label{thm:smooth-periodic-cell-problem}
Let $\gamma\in(0,1)$, $A\in\Gamma_k$, and $h\in C^{2,\gamma}(\mathbb T^n)$ with $h>0$ and $\int_{\mathbb T^n}h=1$. Then \eqref{eq:parametric-corrector} has a unique solution $v_{A,h}\in C^{4,\gamma}(\mathbb T^n)$ satisfying $A+D_y^2v_{A,h}\in\Gamma_k$. If $A\in\mathcal K_0\Subset\Gamma_k$, $h\geq\lambda_0>0$, and $\|h\|_{C^{2,\gamma}}\leq M_0$, then
\[
\|v_{A,h}\|_{C^{4,\gamma}(\mathbb T^n)}\leq C,\qquad A+D_y^2v_{A,h}(y)\in\mathcal K_1\qquad\forall\,y\in\mathbb T^n,
\]
where $\mathcal K_1\Subset\Gamma_k$ and $C$ depend only on $n,k,\gamma,\mathcal K_0,\lambda_0,M_0$. If $h\in C^\infty(\mathbb T^n)$, then $v_{A,h}\in C^\infty(\mathbb T^n)$.
\end{Theorem}

\begin{proof}
For $k=1$, the equation is the periodic Poisson equation $\Delta v=(\operatorname{tr}A)(h-1)$, and the assertions follow from the zero-mean condition and Schauder estimates. Assume henceforth that $k\geq2$. We first derive a priori estimates for an admissible periodic solution $v\in C^{4,\gamma}$ with zero mean. Set
\[
H:=A+D_y^2v,\qquad F:=\sigma_k^{1/k},\qquad g:=(\sigma_k(A)h)^{1/k},\qquad L:=F^{ij}(H)\partial_{y_i y_j}.
\]
The hypotheses give $c\leq g\leq C$ and $\|g\|_{C^2}\leq C$. Choose $\delta>0$, depending only on $\mathcal K_0$, such that $A-\delta I\in\Gamma_k$ for every $A\in\mathcal K_0$. Since \(F(tH)=tF(H)\) for \(t>0\), differentiation at \(t=1\) gives
\[
F^{ij}(H)H_{ij}=F(H).
\]
By the concavity of $F=\sigma_k^{1/k}$, the tangent plane at \(H\) of $F$ lies above the graph of \(F\). Therefore,
\[
\begin{aligned}
F(A-\delta I)&\leq F(H)+F^{ij}(H)(A-\delta I-H)_{ij}=F^{ij}(H)(A-\delta I)_{ij}.
\end{aligned}
\]
Consequently,
\[
F^{ij}(H)A_{ij}=F^{ij}(H)(A-\delta I)_{ij}+\delta\sum_iF^{ii}(H)\geq F(A-\delta I)+\delta\sum_i F^{ii}(H)\geq\delta\sum_i F^{ii}(H).
\]
Differentiating $F(H)=g$ twice and using concavity again gives
\[
L(\Delta v)\geq\Delta g,\qquad Lv=F^{ij}(H)(H_{ij}-A_{ij})=g-F^{ij}(H)A_{ij}.
\]
For $a>0$, let $y_0$ be a maximum point of $\Delta v-av$ on $\mathbb T^n$. At $y_0$,
\[
0\geq L(\Delta v-av)\geq\Delta g-ag+a\delta\sum_iF^{ii}(H),
\]
so $\sum_iF^{ii}(H(y_0))\leq C_a$. Since
\[
\sum_iF^{ii}(H)=\frac{n-k+1}{k}\,\frac{\sigma_{k-1}(H)}{\sigma_k(H)^{(k-1)/k}},
\]
we obtain $\sigma_{k-1}(H(y_0))\leq C_a$. The Newton inequalities imply
\begin{equation*}
\sigma_1(H)\leq C_{n,k}\frac{\sigma_{k-1}(H)^{k-1}}{\sigma_k(H)^{k-2}}.
\end{equation*}
Indeed, the normalized symmetric functions $s_j:=\sigma_j(H)/\binom nj$ are positive for $1\leq j\leq k$ and log-concave, so $s_{k-1}^{k-1}\geq s_1s_k^{k-2}$. It follows that $\operatorname{tr}H(y_0)\leq C_a$. Set $S:=\sup_{\mathbb T^n}\operatorname{tr}H$. By the choice of $y_0$,
\[
S\leq C_a+a\operatorname{osc}_{\mathbb T^n}v.
\]
The Harnack inequality for the Poisson equation and $\operatorname{tr}H>0$ give
\[
\operatorname{osc}_{\mathbb T^n}v\leq C_n\|\Delta v\|_{L^\infty(\mathbb T^n)}\leq C_n(S+\operatorname{tr}A).
\]
Choose $a>0$ so that $aC_n\leq1/2$. Absorbing $aC_nS$, we obtain $S\leq C$. Since $H\in\Gamma_k\subset\Gamma_2$,
\[
|H|^2=(\operatorname{tr}H)^2-2\sigma_2(H)\leq(\operatorname{tr}H)^2\le C,
\]
where $|H|$ is the Frobenius norm. Thus $\|D_y^2v\|_{L^\infty}\leq C$. Then the periodicity and the zero-mean normalization give $\|v\|_{C^2}\leq C$.

Since $|H|\leq C$, $H\in\Gamma_k$, and $\sigma_k(H)\geq c>0$, the Newton--Maclaurin inequalities give a positive lower bound for every $\sigma_j(H)$, $1\leq j\leq k$. Consequently $H$ lies in a fixed compact set $\mathcal K_1\Subset\Gamma_k$, and $T_{k-1}(H)$ is uniformly positive definite.

The operator $\sigma_k^{1/k}$ is concave on $\Gamma_k$ and uniformly elliptic on a neighbourhood of $\mathcal K_1$. The Evans--Krylov estimate gives a uniform $C^{2,\beta}$ bound for some $\beta>0$. Differentiating the equation in $y$ and applying Schauder estimates, we then obtain a $C^{3,\beta'}$ bound for some $\beta'>0$. The linearized coefficients are therefore Lipschitz and have uniform $C^{0,\gamma}$ bounds. A further application of Schauder estimates to the first and second differentiated equations gives the stated $C^{4,\gamma}$ bound; see \cite{gilbarg2001elliptic}*{Chapter~17}. The same argument gives higher regularity for smoother data.

To prove existence, we apply the method of continuity to
\[
\sigma_k(A+D_y^2v_t)=\sigma_k(A)h_t,\qquad h_t:=1-t+th,\qquad \int_{\mathbb T^n}v_t=0,\qquad 0\leq t\leq1.
\]
At $t=0$, take $v_0=0$. The preceding estimates are uniform in $t$ and give closedness. For openness, \eqref{eq:periodic-cell-compatibility} shows that the nonlinear equation takes values in the zero-mean subspace. The linearization $T_{k-1}^{ij}(A+D_y^2v_t)\partial_{y_i y_j}$ is uniformly elliptic and, by \eqref{eq:newton-divergence}, self-adjoint. Its kernel consists of constants, so the Fredholm alternative and Schauder estimates give an isomorphism from the zero-mean $C^{4,\gamma}$ space to the zero-mean $C^{2,\gamma}$ space. The implicit function theorem therefore gives openness. Finally, the difference of two admissible solutions satisfies a uniformly elliptic linear equation, obtained by integrating the linearization along the segment joining their Hessians in the convex cone $\Gamma_k$. The strong maximum principle and the zero-mean normalization give uniqueness.
\end{proof}

We next study the dependence of $v_{A,h}$ on $(A,h)$. For $0<\gamma<1$ and $\ell\in\mathbb N_0$, define
\[
C_0^{\ell,\gamma}(\mathbb T^n):=\left\{r\in C^{\ell,\gamma}(\mathbb T^n):\int_{\mathbb T^n}r(y)\,dy=0\right\},
\]
and write
\[
X_\gamma:=\Gamma_k\times\left\{h\in C^{2,\gamma}(\mathbb T^n):h>0,\ \int_{\mathbb T^n}h(y)\,dy=1\right\}.
\]
Parameter directions are taken in $\mathcal S^{n\times n}\times C_0^{2,\gamma}(\mathbb T^n)$, with the product norm.

\begin{Lemma}
\label{lem:smooth-cell-dependence}
The corrector map $(A,h)\mapsto v_{A,h}$ is $C^3$ from $X_\gamma$ to $C_0^{2,\gamma}(\mathbb T^n)$ in the Fr\'echet sense. If $A\in\mathcal K_0\Subset\Gamma_k$, $h\geq\lambda_0>0$, and $\|h\|_{C^{2,\gamma}}\leq M_0$, then
\[
\|v_{A,h}\|_{C^{2,\gamma}(\mathbb T^n)}+\sum_{j=1}^3\|D_{(A,h)}^jv_{A,h}\|_{C^{2,\gamma}(\mathbb T^n)}\leq C.
\]
Parameter derivatives are measured in the induced operator norms, and $C$ depends only on $n,k,\gamma,\mathcal K_0,\lambda_0,M_0$.
\end{Lemma}

\begin{proof}
Theorem~\ref{thm:smooth-periodic-cell-problem} gives the uniform corrector bound and places $H=A+D_y^2v_{A,h}$ in a fixed compact subset of $\Gamma_k$. The linearized operator
\[
L_{A,h}\eta=T_{k-1}^{ij}(H)\eta_{ij}
\]
is uniformly elliptic. Since its coefficient matrix is divergence-free, it is also self-adjoint. By the strong maximum principle, its kernel on $C^{2,\gamma}$ consists exactly of constants. The Fredholm alternative and Schauder estimates therefore give an isomorphism
\[
L_{A,h}:C_0^{2,\gamma}(\mathbb T^n)\longrightarrow C_0^{0,\gamma}(\mathbb T^n)
\]
with a uniformly bounded inverse. The uniform inverse bound follows from the energy estimate for zero-mean periodic functions and the periodic Schauder estimate.

For $\zeta\in C_0^{2,\gamma}$, \eqref{eq:periodic-cell-compatibility} shows that
\[
\sigma_k(A+D_y^2\zeta)-\sigma_k(A)h\in C_0^{0,\gamma}.
\]
Applying the implicit function theorem at $(A,h,v_{A,h})$, we obtain smooth dependence of $v_{A,h}$ on $(A,h)$. Uniqueness of the normalized cell solution makes these local branches agree. Differentiating \eqref{eq:parametric-corrector} up to three times gives equations with the same operator $L_{A,h}$ and right-hand sides involving parameter directions and parameter derivatives of lower order. We bound these derivatives successively using the uniform inverse bound and the H\"older product estimate.
\end{proof}

\begin{proof}[Proof of Theorem~\ref{thm:main-result}]
We divide the proof into three steps.

\smallskip

\noindent \textbf{Step 1. Reduction and mollification.}
It suffices to consider sufficiently small $\varepsilon$. Indeed, choose a ball $B_R(x_0)$ containing $\Omega$, with $R\leq2\operatorname{diam}\Omega$, and set
\[
b_+=\frac12\left(\frac{\Lambda}{\binom nk}\right)^{1/k},\qquad b_-=\frac12\left(\frac{\lambda}{\binom nk}\right)^{1/k}.
\]
The two quadratic comparison functions give
\[
\min_{\partial\Omega}\varphi+b_+(|x-x_0|^2-R^2)\leq u^\varepsilon(x)\leq\max_{\partial\Omega}\varphi+b_-|x-x_0|^2.
\]
Their $k$-Hessians are $\Lambda$ and $\lambda$, respectively, and their boundary values have the required order. Thus $u^\varepsilon$ is uniformly bounded. Increasing $C$ gives the desired estimate for $\varepsilon$ bounded away from zero.

Set
\[
A(x):=D^2u(x),\qquad h(x,y):=\frac{f(x,y)}{\overline f(x)},\qquad Z(x):=(A(x),h(x,\cdot))\in X_\gamma.
\]
The bounds on $f$ give $\lambda\leq\overline f\leq\Lambda$, $c\leq h\leq C$, and $\int_{\mathbb T^n}h(x,y)\,dy=1$. Since $u$ is $k$-admissible, $\sigma_k(A)=\overline f\geq\lambda$, and $|A|\leq\|u\|_{C^2(\overline\Omega)}$, the Newton--Maclaurin inequalities give
\begin{equation*}
A(x)\in\mathcal K_0\Subset\Gamma_k\qquad\forall\,x\in\overline\Omega.
\end{equation*}
By the convexity of $\Gamma_k$, we may choose $\mathcal K_0$ convex with $A(\overline\Omega)$ in its interior, depending only on $n$, $k$, $\lambda$, and $\|u\|_{C^2(\overline\Omega)}$. Let $E_\gamma:=\mathcal S^{n\times n}\times C^{2,\gamma}(\mathbb T^n)$, equipped with the product norm $\|(B,r)\|_{E_\gamma}:=|B|+\|r\|_{C^{2,\gamma}(\mathbb T^n)}$. Specifically,
\[
\|Z\|_{L^\infty(\Omega;E_\gamma)}=\sup_{x\in\Omega}\left(|A(x)|+\|h(x,\cdot)\|_{C^{2,\gamma}(\mathbb T^n)}\right).
\]
For derivatives of parameter fields, we apply the $E_\gamma$ norm componentwise. Since $\overline f\geq\lambda$, the normalization $h=f/\overline f$ preserves the $C^{m,\alpha}$ regularity in $x$, with values in $C^{2,\gamma}(\mathbb T^n)$, for $m=0,1$.

Using a bounded Whitney extension operator, we extend $A$ and the $C_0^{2,\gamma}(\mathbb T^n)$-valued field $h-1$ in $C^{m,\alpha}$ to a neighborhood of $\overline\Omega$, preserving the normalization of $h$. Shrinking the neighborhood if necessary, we have $A\in\mathcal K_0$ and a uniform positive lower bound for $h$. We then mollify in $x$ with a standard nonnegative symmetric mollifier $\eta_\varepsilon(x)$ at scale $\varepsilon$ and set
\[
A_\varepsilon:=A*\eta_\varepsilon,\qquad h_\varepsilon:=h*\eta_\varepsilon,\qquad Z_\varepsilon:=(A_\varepsilon,h_\varepsilon)=Z*\eta_\varepsilon.
\]
For sufficiently small $\varepsilon$, the convexity of $\mathcal K_0$ and the nonnegativity of the mollifier give $A_\varepsilon(x)\in\mathcal K_0\Subset\Gamma_k$ and $c\leq h_\varepsilon\leq C$. In particular, $Z_\varepsilon:\overline\Omega\to X_\gamma$. Denote the corresponding cell potential on $\mathbb R^n$ and its periodic part by
\[
w^\varepsilon(x,y):=w_{A_\varepsilon(x),h_\varepsilon(x,\cdot)}(y),\qquad v^\varepsilon(x,y):=v_{A_\varepsilon(x),h_\varepsilon(x,\cdot)}(y).
\]
By Theorem~\ref{thm:smooth-periodic-cell-problem}, the matrices $D_y^2w^\varepsilon=A_\varepsilon+D_y^2v^\varepsilon$ take values in a fixed compact set $\mathcal K_1\Subset\Gamma_k$. Fix a compact neighbourhood $\mathcal K_2\Subset\Gamma_k$ of $\mathcal K_1$. These sets and the estimates below depend only on the quantities in Theorem~\ref{thm:main-result}.

\smallskip

\noindent \textbf{Step 2. The case $m=0$.}
Since $Z$ is $C^{0,\alpha}$ in $x$, standard mollification estimates give
\begin{equation}
\label{eq:mollified-hessian}
\|Z_\varepsilon-Z\|_{L^\infty(\Omega;E_\gamma)}+\varepsilon\|DZ_\varepsilon\|_{L^\infty(\Omega;E_\gamma)}+\varepsilon^2\|D^2Z_\varepsilon\|_{L^\infty(\Omega;E_\gamma)}\leq C\varepsilon^\alpha.
\end{equation}

Set
\[
V^\varepsilon(x):=u(x)+\varepsilon^2v^\varepsilon(x,x/\varepsilon).
\]
Differentiating $V^\varepsilon$ twice and writing $y=x/\varepsilon$, we obtain
\[
D_{ij}V^\varepsilon=A_{ij}+(v^\varepsilon)_{y_i y_j}+\varepsilon\bigl((v^\varepsilon)_{x_i y_j}+(v^\varepsilon)_{x_j y_i}\bigr)+\varepsilon^2(v^\varepsilon)_{x_i x_j}.
\]
Since $D_y^2w^\varepsilon=A_\varepsilon+D_y^2v^\varepsilon$, we have
\[
D^2V^\varepsilon(x)=D_y^2w^\varepsilon(x,x/\varepsilon)+E_\varepsilon(x),
\]
with
\[
\begin{aligned}
\|E_\varepsilon\|_{L^\infty(\Omega)}&\leq\|A-A_\varepsilon\|_{L^\infty(\Omega)}+C\varepsilon\|D_xD_yv^\varepsilon\|_{L^\infty(\Omega\times\mathbb T^n)}+\varepsilon^2\|D_x^2v^\varepsilon\|_{L^\infty(\Omega\times\mathbb T^n)}
\leq C\varepsilon^\alpha,
\end{aligned}
\]
where we used \eqref{eq:mollified-hessian} and the estimates
\[
\begin{aligned}
\|D_xD_yv^\varepsilon\|_{L^\infty(\Omega\times\mathbb T^n)}&\leq\|D_xv^\varepsilon\|_{L^\infty(\Omega;C^{2,\gamma}(\mathbb T^n))}\leq C\|DZ_\varepsilon\|_{L^\infty(\Omega;E_\gamma)},\\
\|D_x^2v^\varepsilon\|_{L^\infty(\Omega\times\mathbb T^n)}&\leq C\left(\|D^2Z_\varepsilon\|_{L^\infty(\Omega;E_\gamma)}+\|DZ_\varepsilon\|_{L^\infty(\Omega;E_\gamma)}^{\,2}\right),
\end{aligned}
\]
which follow from Lemma~\ref{lem:smooth-cell-dependence} and the chain rule.

Thus, for sufficiently small $\varepsilon$,
\[
\|V^\varepsilon-u\|_{L^\infty(\overline\Omega)}\leq C\varepsilon^2,\qquad D^2V^\varepsilon(x)\in\mathcal K_2\qquad\forall\,x\in\Omega.
\]
By the cell equation \eqref{eq:parametric-corrector},
\[
\sigma_k(D_y^2w^\varepsilon(x,x/\varepsilon))=\sigma_k(A_\varepsilon(x))h_\varepsilon(x,x/\varepsilon).
\]
Recalling that $f=\sigma_k(A)h$, we use the bound on $E_\varepsilon$, \eqref{eq:mollified-hessian}, and the local Lipschitz continuity of the $k$-Hessian polynomial to obtain
\[
\left|\sigma_k(D^2V^\varepsilon)-f(x,x/\varepsilon)\right|\leq C_0\varepsilon^\alpha.
\]
Since $\varepsilon^2\leq\varepsilon^\alpha$ for $0<\varepsilon\leq1$, Proposition~\ref{prop:approximate-solution-principle}, applied with $\beta=\alpha$, proves the case $m=0$.

\smallskip

\noindent  \textbf{Step 3. The case $m=1$.}
For each $x\in\overline\Omega$, let $\psi^\varepsilon(x,\cdot)$ be the zero-mean solution of
\begin{equation}
\label{eq:third-corrector}
(W^\varepsilon)^{ij}(\psi^\varepsilon)_{y_i y_j}=-2(W^\varepsilon)^{ij}(v^\varepsilon)_{x_i y_j}\quad\text{in }\mathbb T^n,\qquad\int_{\mathbb T^n}\psi^\varepsilon(x,y)\,dy=0,
\end{equation}
where $W^\varepsilon(x,y):=T_{k-1}(D_y^2w^\varepsilon(x,y))$. Since $W^\varepsilon$ is divergence-free in $y$,
\[
\int_{\mathbb T^n}(W^\varepsilon)^{ij}(v^\varepsilon)_{x_i y_j}\,dy=0.
\]
The uniform ellipticity of $W^\varepsilon$ and the Fredholm alternative therefore give a unique zero-mean solution.

We next estimate the correctors and their derivatives in $x$. The $C^{1,\alpha}$ regularity of $Z$ gives
\begin{equation}
\label{eq:c1alpha-mollification}
\begin{alignedat}{2}
\|Z_\varepsilon-Z\|_{L^\infty(\Omega;E_\gamma)}&\leq C\varepsilon^{1+\alpha},\qquad&\|DZ_\varepsilon-DZ\|_{L^\infty(\Omega;E_\gamma)}&\leq C\varepsilon^\alpha,\\
\|D^2Z_\varepsilon\|_{L^\infty(\Omega;E_\gamma)}&\leq C\varepsilon^{\alpha-1},\qquad&\|D^3Z_\varepsilon\|_{L^\infty(\Omega;E_\gamma)}&\leq C\varepsilon^{\alpha-2}.
\end{alignedat}
\end{equation}
The first estimate uses the symmetry of $\eta_\varepsilon$. The other three follow from standard mollification estimates for $DZ\in C^{0,\alpha}$. Lemma~\ref{lem:smooth-cell-dependence} and the chain rule give
\[
\begin{alignedat}{2}
\|v^\varepsilon\|_{L^\infty(\Omega;C^{2,\gamma}(\mathbb T^n))}&\leq C,\qquad&\|D_xv^\varepsilon\|_{L^\infty(\Omega;C^{2,\gamma}(\mathbb T^n))}&\leq C,\\
\|D_x^2v^\varepsilon\|_{L^\infty(\Omega;C^{2,\gamma}(\mathbb T^n))}&\leq C\varepsilon^{\alpha-1},\qquad&\|D_x^3v^\varepsilon\|_{L^\infty(\Omega;C^{2,\gamma}(\mathbb T^n))}&\leq C\varepsilon^{\alpha-2}.
\end{alignedat}
\]
Since $W^\varepsilon=T_{k-1}(D_y^2w^\varepsilon)$, these estimates also give
\[
\|W^\varepsilon\|_{L^\infty(\Omega;C^{0,\gamma}(\mathbb T^n))}+\|D_xW^\varepsilon\|_{L^\infty(\Omega;C^{0,\gamma}(\mathbb T^n))}+\varepsilon^{1-\alpha}\|D_x^2W^\varepsilon\|_{L^\infty(\Omega;C^{0,\gamma}(\mathbb T^n))}\leq C.
\]
Applying Schauder estimates to \eqref{eq:third-corrector} and to the equations obtained by differentiating it once and twice with respect to $x$, we obtain
\[
\|\psi^\varepsilon\|_{L^\infty(\Omega;C^{2,\gamma}(\mathbb T^n))}+\varepsilon^{1-\alpha}\|D_x\psi^\varepsilon\|_{L^\infty(\Omega;C^{2,\gamma}(\mathbb T^n))}+\varepsilon^{2-\alpha}\|D_x^2\psi^\varepsilon\|_{L^\infty(\Omega;C^{2,\gamma}(\mathbb T^n))}\leq C.
\]

Set
\[
\widetilde V^\varepsilon(x):=u(x)+\varepsilon^2v^\varepsilon(x,x/\varepsilon)+\varepsilon^3\psi^\varepsilon(x,x/\varepsilon).
\]
Differentiating twice, we obtain
\begin{equation}
\label{eq:third-order-hessian}
D^2\widetilde V^\varepsilon(x)=D_y^2w^\varepsilon(x,x/\varepsilon)+\varepsilon N_\varepsilon(x,x/\varepsilon)+R_\varepsilon(x,x/\varepsilon),
\end{equation}
where
\[
(N_\varepsilon)_{ij}(x,y):=(v^\varepsilon)_{x_i y_j}+(v^\varepsilon)_{x_j y_i}+(\psi^\varepsilon)_{y_i y_j}
\]
and
\[
(R_\varepsilon)_{ij}(x,y):=(A-A_\varepsilon)_{ij}+\varepsilon^2(v^\varepsilon)_{x_i x_j}+\varepsilon^2\bigl((\psi^\varepsilon)_{x_i y_j}+(\psi^\varepsilon)_{x_j y_i}\bigr)+\varepsilon^3(\psi^\varepsilon)_{x_i x_j}.
\]
The preceding derivative bounds give
\[
\|R_\varepsilon\|_{L^\infty(\Omega\times\mathbb T^n)}\leq \|A-A_\varepsilon\|_{L^\infty(\Omega)}+C\bigl(\varepsilon^2\varepsilon^{\alpha-1}+\varepsilon^3\varepsilon^{\alpha-2}\bigr)\leq C\varepsilon^{1+\alpha}.
\]
Moreover, $N_\varepsilon$ and $D_y^2w^\varepsilon$ are uniformly bounded, and $D_y^2w^\varepsilon$ takes values in $\mathcal K_1$. By \eqref{eq:third-order-hessian} and the bound on $R_\varepsilon$, $D^2\widetilde V^\varepsilon(x)\in\mathcal K_2$ for all $x\in\Omega$ and sufficiently small $\varepsilon$. By \eqref{eq:third-corrector},
\[
(W^\varepsilon)^{ij}(N_\varepsilon)_{ij}=0\qquad\text{in }\Omega\times\mathbb T^n.
\]
Thus the term of order $\varepsilon$ in the $k$-Hessian polynomial expansion vanishes. Using the cell equation \eqref{eq:parametric-corrector}, \eqref{eq:c1alpha-mollification}, and the bound on $R_\varepsilon$, we obtain
\[
\begin{aligned}
\left|\sigma_k(D^2\widetilde V^\varepsilon)-f(x,x/\varepsilon)\right|\leq C\bigl(\varepsilon^2+\varepsilon^{1+\alpha}\bigr)+C\|\sigma_k(A_\varepsilon)h_\varepsilon-\sigma_k(A)h\|_{L^\infty(\Omega\times\mathbb T^n)}\leq C\varepsilon^{1+\alpha}.
\end{aligned}
\]
Finally, the uniform bounds for $v^\varepsilon$ and $\psi^\varepsilon$ give
\[
\|\widetilde V^\varepsilon-u\|_{L^\infty(\overline\Omega)}\leq C\varepsilon^2.
\]
Then Proposition~\ref{prop:approximate-solution-principle} applies with $\beta=1+\alpha$ and completes the proof.
\end{proof}

\textbf{AI usage statement.}  GPT-6 Astra was used for language editing.

\bigskip

\noindent T. Jin

\noindent Department of Mathematics, The Hong Kong University of Science and Technology\\
Clear Water Bay, Kowloon, Hong Kong\\
Email: \textsf{tianlingjin@ust.hk}

\bigskip

\noindent Y.Y. Li

\noindent Department of Mathematics, Rutgers University\\ 
110 Frelinghuysen Road, Piscataway, NJ 08854-8019, USA\\
Email: \textsf{yyli@math.rutgers.edu}

\bigskip

\noindent H.V. Tran

\noindent Department of Mathematics, University of Wisconsin at Madison,\\
Van Vleck Hall, 480 Lincoln Drive, Madison, WI 53706, USA\\
Email: \textsf{hung@math.wisc.edu}

\bigskip

\noindent X. Tu

\noindent Department of Mathematics, The Hong Kong University of Science and Technology\\
Clear Water Bay, Kowloon, Hong Kong\\[1mm]
Email:  \textsf{maxstu@ust.hk}; \textsf{tuxushan28561m@gmail.com}

\end{document}